\documentclass[12pt]{article}
\newtheorem{lemma}{Lemma}[section]
\newtheorem{theorem}[lemma]{Theorem}
\newtheorem{proposition}[lemma]{Proposition}

\newtheorem{definition}[lemma]{Definition}
\newenvironment{proof}{{\bf Proof}}{{\hfill $ \Box $}\vskip 4mm}
\newenvironment{remark}{\addtocounter{lemma}{1}
{\bf Remark \thelemma}}{{\hfill}\vskip 4mm}
\newenvironment{remarks}{\addtocounter{lemma}{1}
{\bf Remarks \thelemma}}{{\hfill}\vskip 4mm}

\newenvironment{examples}{\addtocounter{lemma}{1}
{\bf Examples \thelemma}}{{\hfill}\vskip 4mm}
\newcommand{\nc}{\newcommand}
\nc{\rnc}{\renewcommand}
\nc{\nt}{\newtheorem}

        {~\hfill~\fbox{}\end{trivlist}}

\nc{\thlabel}[1]{\label{theo:#1}}
\nc{\thref}[1]{Theorem~\ref{theo:#1}}
\nc{\selabel}[1]{\label{sect:#1}}
\nc{\seref}[1]{Section~\ref{sect:#1}}
\nc{\lelabel}[1]{\label{lemm:#1}}
\nc{\leref}[1]{Lemma~\ref{lemm:#1}}
\nc{\prlabel}[1]{\label{prop:#1}}
\nc{\prref}[1]{Proposition~\ref{prop:#1}}
\nc{\colabel}[1]{\label{coro:#1}}
\nc{\coref}[1]{Corollary~\ref{coro:#1}}
\nc{\exlabel}[1]{\label{exam:#1}}
\nc{\exref}[1]{Example~\ref{exam:#1}}
\nc{\delabel}[1]{\label{defi:#1}}
\nc{\deref}[1]{Definition~\ref{defi:#1}}
\nc{\eqlabel}[1]{\label{equation:#1}}
\nc{\eqref}[1]{(\ref{equation:#1})}
\nc{\csm}{\mbox{$\triangleright\!\!\!<$}}
\nc{\smc}{\mbox{$>\!\!\!\triangleleft$}}
\nc{\trr}{\triangleright}
\providecommand{\operatorname}[1]{\mathrm{#1}\,}
\nc{\Hom}{\operatorname{Hom}}
\nc{\Mor}{\operatorname{Mor}}
\nc{\Aut}{\operatorname{Aut}}
\nc{\Ann}{\operatorname{Ann}}
\nc{\Ker}{\operatorname{Ker}}
\nc{\Trace}{\operatorname{Trace}}
\nc{\Char}{\operatorname{Char}}
\nc{\Mod}{\operatorname{Mod}}
\nc{\End}{\operatorname{End}}
\nc{\Spec}{\operatorname{Spec}}
\nc{\Span}{\operatorname{Span}}
\nc{\sgn}{\operatorname{sgn}}
\nc{\Id}{\operatorname{Id}}
\nc{\Com}{\operatorname{Com}}

\def\Box{\mbox{$\sqcap\!\!\!\!\sqcup$}}

\nc{\dht}{\mbox{$\rightharpoonup\hspace{-2ex}\rightharpoonup$}}
\nc{\dhtb}{\mbox{$\leftharpoonup\hspace{-2ex}\leftharpoonup$}}
\nc{\nd}{\mbox{$\not|$}} 
\providecommand{\text}[1]{\mbox{{\textrm #1}}}

\nc{\nci}{\mbox{$\not\subseteq$}}
\nc{\scontainin}{\mbox{$\mbox{}\subseteq\hspace{-1.5ex}\raisebox{-.5ex}{$_\prime
$}\hspace*{1.5ex}$}}

\def\ot{\otimes}

\def\doublerightleft#1#2{{\lower.2ex\vbox{
\hbox{${\smash{\mathop{\longrightarrow}\limits^{#1}}}$}\vspace*{-4mm}
\hbox{${\smash{\mathop{\longleftarrow}\limits_{#2}}}$}}}}

\newfont{\bbb}{msbm10 scaled\magstep1}  
\newfont{\bbbsub}{msbm10}                
\newfont{\msam}{msam10 scaled\magstep1}

\begin{document}
\title{The Hopf modules category and the Hopf equation}
\author{G. Militaru
\\University of Bucharest
\\Faculty of Mathematics\\Str. Academiei 14
\\RO-70109 Bucharest 1, Romania
\\e-mail: gmilit@al.math.unibuc.ro}
\date{}
\maketitle
\begin{abstract}
\noindent Let $({\cal A}, \Delta )$ be a Hopf-von Neumann algebra and $R$
be the unitary fundamental operator on ${\cal A}$ defined by Takesaki in
\cite{Ta}: $R(a\ot b)=\Delta(b)(a\ot 1)$. Then
$R^{12}R^{23}=R^{23}R^{13}R^{12}$ (see lemma 4.9 of \cite{Ta}). This
operator $R$ plays a vital role in the theory of duality for von Neumann
algebras (see \cite{Ta} or \cite{BS}). If $V$ is a vector space over
an arbitrary field $k$, we shall study what we have called
the Hopf equation: $R^{12}R^{23}=R^{23}R^{13}R^{12}$ in
$\End_k(V\ot V\ot V)$. Taking $W:=\tau R\tau$, the Hopf equation is
equivalent with the pentagonal equation: $W^{12}W^{13}W^{23}=W^{23}W^{12}$
from the theory of operator algebras (see \cite{BS}), where $W$ are viewed
as map in ${\cal L}(K\ot K)$, for a Hilbert space $K$.
For a bialgebra $H$, we shall prove that the classic category of Hopf modules
${}_H{\cal M}^H$ plays a decisive role in describing all solutions of the
Hopf equation. More precisely, if $H$ is a bialgebra over $k$ and
$(M,\cdot ,\rho) \in {}_H{\cal M}^H$ is an $H$-Hopf module, then
the natural map $R=R_{(M,\cdot ,\rho )}$ is a solution for the
Hopf equation. Conversely, the main result of this paper is a FRT type
theorem: if $M$ is a finite dimensional vector space and
$R\in \End_k(M\ot M)$ is a solution for the Hopf equation, then there
exists a bialgebra $B(R)$ such that
$(M,\cdot ,\rho) \in {}_{B(R)}{\cal M}^{B(R)}$ and $R=R_{(M,\cdot ,\rho )}$.
By applying this result, we construct new examples of noncommutative
and noncocommutative bialgebras which are different from the ones
arising from quantum group theory. In particular, over a field of
characteristic two, an example of five dimensional noncommutative and
noncocommutative bialgebra is given.
\end{abstract}

\section{Introduction}
Let $H$ be a bialgebra over a field $k$. There are two fundamental categories
in the theory of Hopf algebras and quantum groups: ${}_H{\cal M}^H$, the
category of $H$-Hopf modules and ${}_H{\cal YD}^H$, the category of quantum
Yetter-Drinfel'd modules. The objects in these categories are $k$-vector
spaces $M$ which are left $H$-modules $(M,\cdot )$, right $H$-comodules
$(M,\rho )$, such that
the following quite distinct compatibility relations hold:
\begin{equation}\label{H}
\rho(h\cdot m)=\sum h_{(1)}\cdot m_{<0>}\ot h_{(2)}m_{<1>}
\end{equation}
in the case ${}_H{\cal M}^H$, and respectively
\begin{equation}
\sum h_{(1)}\cdot m_{<0>}\ot h_{(2)}m_{<1>}
=\sum (h_{(2)}\cdot m)_{<0>}\ot (h_{(2)}\cdot m)_{<1>}h_{(1)}
\end{equation}
for the Yetter-Drinfel'd categories.

Traditionally, these two categories have been studied for completely
different reasons: the classical category ${}_H{\cal M}^H$
(or immediate generalisations of it: ${}_A{\cal M}^H$, ${}_A{\cal M}(H)^C$)
is involved in the theory of integrals for a Hopf algebra
(see \cite{A}, \cite{Sw} or the more recent \cite{M}), Clifford theory
of representations (\cite{MS}, \cite{Sch1}, \cite{Sch3}, \cite{St})
and Hopf-Galois theory (\cite{M}, \cite{Sch2}, etc.). The cateory
${}_H{\cal YD}^H$, introduced in \cite{Y}, plays an important role in the
quantum Yang-Baxter equation, quantum groups, low dimensional topology and
knot theory (see \cite{K}, \cite{LR}, \cite{R1}, \cite{RT},
or \cite{T}).

However, there are two connections between these categories. The first
one was given by P. Schauenburg in \cite{S1}: it was proven that the
category ${}_H{\cal YD}^H$ is equivalent to the category
${}_H^H{\cal M}^H_H$ of two-sided, two-cosided Hopf modules. The second
was given recentely in \cite{CMZ1}. For $A$ an $H$-comodule algebra
and $C$ an $H$-module coalgebra, Doi (see \cite{D}) and independently
Koppinen (see \cite{Ko}) defined ${}_A{\cal M}(H)^C$, the category of
Doi-Koppinen Hopf modules, whose objects
are left $A$-modules and right $C$-comodules and satisfy a compatibility
relation which generalises (\ref{H}). In \cite{CMZ1} it was proven that
${}_H{\cal YD}^H$ is isomorphic to ${}_H{\cal M}(H^{op}\ot H)^H$, where
$H$ can be viewed as an $H^{op}\ot H$-module (comodule) coalgebra (algebra).
The isomorhism is the identity functor $M\rightarrow M$. We hereby
obtain a strong link between the categories ${}_H{\cal M}^H$ and
${}_H{\cal YD}^H$: both are particular cases of the same general category
${}_A{\cal M}(H)^C$. This led us in \cite{CMZ2}, \cite{CMZ3} to study
the implications of the category ${}_H{\cal YD}^H$ in the classic,
non-quantic part of Hopf algebra theory.
In \cite{CMZ2} we start with the following classic theorem (see \cite{M}):
any finite dimensional Hopf algebra is Frobenius. In the language of
categories, this result is interpreted as follows: the forgetful functor
${}_H{\cal M}^H\rightarrow {}_H{\cal M}$
is Frobenius (i.e., cf. \cite{CMZ2}, by definition has the same left and
right adjoint)
if and only if H is finite dimensional. The next step is easy to take: we
must generalize this result for the forgetful functor
${}_A{\cal M}(H)^C\rightarrow {}_A{\cal M}$
and then apply it
in the case of Yetter-Drinfel'd modules for the forgetful functor
${}_H{\cal YD}^H\rightarrow {}_H{\cal M}$.
We thus obtain the fact that the forgetful functor
${}_H{\cal YD}^H\rightarrow {}_H{\cal M}$ is Frobenius if and
only if H is finite dimensional and unimodular (see theorem 4.2 of
\cite{CMZ2}). The same treatment was applied in \cite{CMZ3} for the
classic Maschke theorem.
One of the major obstacles was to correctly define the notion of integral
for the Doi-Hopf datum $(H,A,C)$, such as to be connected to the classic
integral on a Hopf algebra (corresponding to the case C=A=H), as well as to
the notion of total integral (corresponding to the case C=A)
defined by Doi in \cite{D2}.
This technique can be looked upon as a "quantisation" of the theorems from
the classic theory of Hopf algebras. There are two steps to it: first, we
seek to generalize a result for the category ${}_A{\cal M}(H)^C$,
then to apply it to the particular ${}_H{\cal YD}^H$ case. There is also
another approach to this "quantisation" technique, recently evidenced in
\cite{FMS} for the same Frobenius type theorem. It was first proven, by
generalizing the classic result, that any finite dimensional Hopf
algebras extensions is a $\beta$-Frobenius extension (or a Frobenius
extension of second kind). Then, this theorem was "quantised" to the case
of Hopf algebras extensions in ${}_H{\cal YD}^H$. The result includes
the case of enveloping algebras of Lie coloralgebras.

Beginning with this paper, we shall tackle the reverse problem: we shall
try to involve the category ${}_H{\cal M}^H$ in fields dominated until
now by ${}_H{\cal YD}^H$, i.e. try a "dequantisation".
For the beginning, it is enough to remind that the category
${}_H{\cal YD}^H$ is deeply involved in the quantum Yang-Baxter equation:
\begin{equation}
R^{12}R^{13}R^{23}=R^{23}R^{13}R^{12}
\end{equation}
where $R\in \End_k(M\ot M)$, $M$ being a $k$-vector space.
The starting point of this paper is the following question:

{\sl "Can the category ${}_H{\cal M}^H$ be studied in connection with a
certain non-linear equation?"}

The answer is affirmative and, surprisingly, the equation in which the
category ${}_H{\cal M}^H$ is involved (which we shall call
{\sl Hopf equation}) is very close to the quantum Yang-Baxter equation.
More precisely, it is
\begin{equation}
R^{12}R^{23}=R^{23}R^{13}R^{12}
\end{equation}
The simple way of obtaining it from the quantum Yang-Baxter equation by
just deleting the term $R^{13}$ on the left hand side has nevertheless
unpleasant effects: first of all, if the Yang-Baxter equation is reduced
to the solution of a homogenous system, in the case of the Hopf equation the
system is not homogenous any more; secondly, if $R$ is a  solution
of the Hopf equation, $W:=\tau R\tau $ (or $W:=R^{-1}$, if $R$ is bijective)
is not a solution for the Hopf equation, but for the pentagonal equation:
$$
W^{12}W^{13}W^{23}=W^{23}W^{12}.
$$
An interesting connection between the pentagonal and the quantum
Yang-Baxter equations is given in \cite{VV}.
The pentagonal equation plays a fundamental role in the duality theory
for operator algebras (see \cite{BS} and the references indicated here).
If $H$ is a Hopf algebra, then
$$
R:H\ot H\to H\ot H,\quad R(g\ot h)=\sum h_{(1)}g\ot h_{(2)}
$$
is a bijective solution of the Hopf equation. Moreover, the
comultiplication $\Delta$ can be rebuilt from $R$ by means of
$$
\Delta(h)=R((1\ot h)z),
$$
where $z\in H\ot H$ such that $R(z)=1\ot 1$. This operator was defined
first by Takesaki in \cite{Ta} for a Hopf-von Neumann algebra.
The operator $W:=\tau R\tau$ is called in \cite{Maj 1} the evolution operator for a Hopf algebra and
plays an important role in the description of the Markov transition operator
for the quantum random walks (see \cite{Maj 1} or \cite{Maj 2}).

The Hopf equation can be viewed as a natural generalisation of the
idempotent endomorphisms of a vector space: more precisely, if
$f\in \End_k(M)$, then $f\ot I$ (or $I\ot f$)
is a solution of the Hopf equation if and only if $f^2=f$.
We shall prove that if $(M,\cdot ,\rho) \in {}_H{\cal M}^H$
then the natural map
$$R_{(M,\cdot ,\rho )}(m\ot n)=\sum n_{<1>}\cdot m\ot n_{<0>}$$
is a solution of the Hopf equation. Conversely, the main result of this
paper is a FRT type theorem which shows that in the finite dimensional
case, any solution $R$ of the Hopf equation has this form, i.e. there
exists a bialgebra $B(R)$ such that
$(M,\cdot ,\rho) \in {}_{B(R)}{\cal M}^{B(R)}$ and $R=R_{(M,\cdot ,\rho )}$.
Similarly to \cite{BS}, a solution $R$ of the Hopf equation is called
commutative if $R^{12}R^{13}=R^{13}R^{12}$. In the finite dimensional case,
any commutative solution of the Hopf equation has the form
$R=R_{(M,\cdot ,\rho )}$, where $(M,\cdot ,\rho )$ is a Hopf module over a
commutative bialgebra $\overline{B}(R)$. This result can be viewed as the
algebraic version of the theorem 2.2 of \cite{BS}, which classifies all
multiplicative, unitary and commutative operators which can be defined
on a Hilbert space.
In the last part we shall apply our theorem to constructing new examples
of noncommutative, noncocommutative bialgebras which differ from those
arising from the FRT theorem for the quantum Yang-Baxter equation.
These bialgebras arise from the elementary maps
of plane euclidian geometry: the projections of $k^2$ on the $Ox$ and
$Oy$ coordinate axis.
Surprisingly, over a field of characteristic 2, our FRT type construction
supplies us an example of noncommutative and noncocommutative bialgebra of
dimension 5.

Obviously, substituting the key map $R_{(M,\cdot, \rho)}$ with
$\tau R_{(M,\cdot, \rho)} \tau$ , all results of this paper
remain valid if we replace the Hopf equation with the pentagonal equation.
We have preferred however to work with the Hopf equation, for historical
reasons: this is how the issue has been raised for the first time in lemma
4.9 of \cite{Ta}.

This study was continued in \cite{Mi}, where new classes of bialgebras
arising from the Hopf equation are introduced and analyzed.

\section{Preliminaries}
Throughout this paper, $k$ will be a field.
All vector spaces, algebras, coalgebras and bialgebras that we consider
are over $k$. $\ot$ and $\Hom$ will mean $\ot_k$ and $\Hom_k$.
For a coalgebra $C$, we will use Sweedler's $\Sigma$-notation, that is,
$\Delta(c)=\sum c_{(1)}\ot c_{(2)},~(I\ot\Delta)\Delta(c)=
\sum c_{(1)}\ot c_{(2)}\ot c_{(3)}$, etc. We will
also use  Sweedler's notation for right $C$-comodules:
$\rho_M(m)=\sum m_{<0>}\otimes m_{<1>}$, for any $m\in M$ if
$(M,\rho_M)$ is a right $ C$-comodule. ${\cal M}^C$ will be the
category of right $C$-comodules and $C$-colinear maps and
${}_A{\cal M}$ will be the category of left $A$-modules and
$A$-linear maps, if $A$ is a $k$-algebra.

Recall the following well known lemmas:

\begin{lemma}\label{unu}
Let $M$ be a finite dimensional vector space with $\{m_1,\cdots, m_n \}$
a basis for $M$ and let $C$ be a coalgebra. We define the $k$-linear map
$\rho :M\to M\ot C$, $\rho(m_l)=\sum_{v=1}^{n}m_v\ot c_{vl}$, for all
$l=1,\cdots, n$, where $(c_{vl})_{v,l}$ is a family of elements of $C$.
The following statements are equivalent:
\begin{enumerate}
\item $(M,\rho)$ is a right $C$-comodule.
\item The matrix $(c_{vl})_{v,l}$ is comultiplicative, i.e.
\begin{equation}\label{com}
\Delta(c_{jk})=\sum_{u=1}^{n}c_{ju}\ot c_{uk},
\quad \varepsilon (c_{jk})=\delta_{jk}
\end{equation}
for all $j,k=1,\cdots, n$
\end{enumerate}
\end{lemma}

If we denote $B=(c_{vl})_{v,l}$, then, as usual, the relations
(\ref{com}) can formally be written:
$\Delta (B)=B\ot B$, $\varepsilon (B)=I_n$.

\begin{lemma}\label{doi}
Let $(C,\Delta, \varepsilon)$ be a coalgebra. Then, on the tensor algebra
$(T(C), M, u)$, there exists a unique bialgebra structure
$(T(C), M, u, \overline{\Delta},\overline{\varepsilon})$ such that
$\overline{\Delta}(c)=\Delta (c)$ and
$\overline{\varepsilon}(c)=\varepsilon (c)$ for all $c\in C$.
In addition, the inclusion map $i:C\to T(C)$ is a coalgebra
map.\\
Furthermore, if $M$ is a vector space and $\mu: C\ot M\to M$,
$\mu (c\ot m)=c\cdot m$ is a linear map, then there exists a unique left
$T(C)$-module structure on $M$, $\overline{\mu}:T(C)\ot M\to M$, such that
$\overline{\mu}(c\ot m)=c\cdot m$, for all $c\in C$, $m\in M$.
\end{lemma}

Let $H$ be a bialgebra. Recall that an (left-right) $H$-Hopf module is a
left $H$-module $(M, \cdot)$ which is also a right $H$-comodule $(M,\rho)$
such that
\begin{equation}\label{H1}
\rho(h\cdot m)=\sum h_{(1)}\cdot m_{<0>}\ot h_{(2)}m_{<1>}
\end{equation}
for all $h\in H$, $m\in M$. ${}_H{\cal M}^H$ will be the category of
$H$-Hopf modules and $H$-linear $H$-colinear homomorphisms.

\begin{lemma} \label{trei}
Let $H$ be a bialgebra, $(M,\cdot)$ a left $H$-module
and $(M,\rho)$ a right $H$-comodule. Then the set
$$
\{h\in H\mid \rho(h\cdot m)=\sum h_{(1)}\cdot m_{<0>}\ot h_{(2)}m_{<1>},
\forall m\in M \}
$$
is a subalgebra of $H$.
\end{lemma}

\begin{proof} Straightforward.
\end{proof}

We obtain from this lemma that if a left $H$-module and  right
$H$-comodule $M$ satisfies the condition of compatibility (\ref{H1})
for a set of generators as an algebra of $H$ and for a basis of $M$, then
$M$ is an $H$-Hopf module. If $(M,\cdot)$ is a left $H$-module and
$(M,\rho)$ is a right $H$-comodule, the special map
\begin{equation}\label{sp}
R_{(M,\cdot ,\rho )}:M\ot M\to M\ot M, \quad
R_{(M,\cdot ,\rho )}(m\ot n)=\sum n_{<1>}\cdot m\ot n_{<0>}
\end{equation}
will play an important role in the present paper. It is useful to point out
the following lemma. The proof is left to the reader.

\begin{lemma}\label{patru}
Let $H$ be a bialgebra, $(M,\cdot)$ a left $H$-module
and $(M,\rho)$ a right $H$-comodule. If $I$ is a biideal of $H$ such that
$I\cdot M=0$, then, with the natural structures, $(M,\cdot^{\prime})$ is a
left $H/I$-module, $(M,\rho^{\prime})$ a right $H/I$-comodule and
$R_{(M,\cdot^{\prime},\rho^{\prime})}=R_{(M,\cdot ,\rho )}$
\end{lemma}

For a vector space $V$, $\tau :V\otimes V\to V\otimes V$
will denote the switch map, that is, $\tau (v\otimes w)=w\otimes v$
for all $v,w \in V$. If $R:V\ot V\to V\ot V$ is a linear map
we denote by $R^{12}$, $R^{13}$, $R^{23}$ the maps of $\End_k(V\ot V\ot V)$
given by
$$
R^{12}=R\ot I, \quad R^{23}=I\ot R,\quad
R^{13}=(I\ot \tau)(R\ot I)(I\ot \tau).
$$
Using the notation $R(u\ot v)=\sum u_1\ot v_1$ then
$$
R^{12}(u\ot v\ot w)=\sum u_1\ot v_1\ot w_0
$$
where the subscript (0) means that $w$ is not affected by the application
of $R^{12}$.\\
Let $H$ be a bialgebra and $(M, \cdot)$ a left $H$-module which is also a
right $H$-comodule $(M,\rho)$. Recall that $(M,\cdot,\rho)$ is a
Yetter-Drinfel'd module if the following compatibility relation holds:
$$
\sum h_{(1)}\cdot m_{<0>}\ot h_{(2)}m_{<1>}
=\sum (h_{(2)}\cdot m)_{<0>}\ot (h_{(2)}\cdot m)_{<1>}h_{(1)}
$$
for all $h\in H$, $m\in M$. ${}_H{\cal YD}^H$ will be the category of
Yetter-Drinfel'd modules and $H$-linear $H$-colinear homomorphism.
If $(M,\cdot ,\rho)$ is a Yetter-Drinfel'd module then the special map
$R=R_{(M,\cdot ,\rho )}$ given by the equation (\ref{sp}) is a solution
of the quantum Yang-Baxter equation
$$
R^{12}R^{13}R^{23}=R^{23}R^{13}R^{12}.
$$
If $M$ is a finite dimensional vector space and $R$ is a solution of the
quantum Yang-Baxter equation, then there exists a bialgebra $A(R)$
such that $(M,\cdot ,\rho) \in {}_{A(R)}{\cal M}^{A(R)}$ and
$R=R_{(M,\cdot ,\rho )}$ (see \cite{R1}). For a further study of the
Yetter-Drinfel'd category we refer to \cite{LR}, \cite{R1}, \cite{RT},
\cite{Y}, or to the more recent \cite{CMZ1}, \cite{CMZ2}, \cite{CMZ3},
\cite{FMS}.

\section{The Hopf equation}
We will start with the following

\begin{definition}
Let $V$ be a vector space and $R\in \End_k(V\ot V)$.
\begin{enumerate}
\item We shall say that $R$ is a solution for the Hopf equation if
\begin{equation}\label{Heq}
R^{23}R^{13}R^{12}=R^{12}R^{23}
\end{equation}
\item We shall say that $R$ is a solution for the pentagonal equation if
\begin{equation}\label{iHeq}
R^{12}R^{13}R^{23}=R^{23}R^{12}
\end{equation}
\end{enumerate}
\end{definition}

\begin{remarks}
1. The Hopf equation is obtained from the quantum Yang-Baxter equation
$$R^{23}R^{13}R^{12}=R^{12}R^{13}R^{23}$$
by deleting the midle term from the right hand side.

2. Let $\{m_i \}_{i\in I}$ be a basis of the vector space $V$. Then an
endomorphism $R$ of $V\ot V$ is given by a family of scalars
$(x_{ij}^{kl})_{i,j,k,l\in I}$ of $k$ such that
\begin{equation}\label{sc}
R(m_v\ot m_u)=\sum _{i,j}x_{uv}^{ji}m_i\ot m_j
\end{equation}
for all $v,u\in I$. A direct computation shows us that $R$ is a solution of
the Hopf equation if and only if $(x_{ij}^{kl})_{i,j,k,l\in I}$ is a
solution of the nonlinear equation
\begin{equation}\label{neq}
\sum_{\alpha ,\beta ,\gamma}
x_{\gamma \alpha}^{jk}x_{w\beta}^{\gamma l}x_{uv}^{\alpha \beta}=
\sum_{i} x_{wu}^{ji}x_{iv}^{kl}
\end{equation}
for all $j, k, l, u, v, w\in I$. It follows that solving the system
(\ref{neq}) is really a non-trivial problem.

3. Using the notation $R(x\ot y)=\sum x_1\ot y_1$, for $x, y\in V$, then
$R$ is a solution of the Hopf equation if and only if
$$
\sum x_{110}\ot y_{101}\ot z_{011}=
\sum x_{01}\ot y_{11}\ot z_{10}
$$
for all $x$, $y$, $z\in V$.

4. Suppose that $R\in \End_k(V\ot V)$ is bijective. Then, $R$ is a solution
of the Hopf equation if and only if $R^{-1}$ is a solution of the
pentagonal equation.

5. Let $A$ be an algebra and $R\in A\ot A$ be an invertible element such
that the Hopf equation  $R^{23}R^{13}R^{12}=R^{12}R^{23}$ holds in
$A\ot A\ot A$. Then, the comultiplication
$$\Delta :A\to A\ot A,\quad \Delta(a):=R(1\ot a)R^{-1}$$
for all $a\in A$ is coassociative and an algebra map.

Indeed,
$$
(I\ot \Delta)\Delta (a)=R^{23}R^{13}(1\ot 1\ot a)(R^{23}R^{13})^{-1}
$$
and
$$
(\Delta\ot I)\Delta (a)=R^{12}R^{23}(1\ot 1\ot a)(R^{12}R^{23})^{-1}
$$
Let $W:=(R^{13})^{-1}(R^{23})^{-1}R^{12}R^{23}$. Then $\Delta$ is
coassociative if and only if
\begin{equation}\label{jaz}
(1\ot 1\ot a)W=W(1\ot 1\ot a)
\end{equation}
for all $a\in A$. But, as $R$ satisfies the Hopf equation, we have that
$W=R^{12}$, i.e. equation (\ref{jaz}) holds.
\end{remarks}

In the next proposition we shall evidence a few equations which are
equivalent to the Hopf equation.

\begin{proposition}
Let $V$ be a vector space and $R\in \End_k(V\ot V)$. The following
statements are equivalent:
\begin{enumerate}
\item $R$ is a solution of the Hopf equation.
\item $T:=\tau R$ is a solution of the equation:
$T^{12}T^{23}T^{12}=T^{23}\tau^{12}T^{23}$.
\item $T:=R\tau$ is a solution of the equation:
$T^{23}T^{12}T^{23}=T^{12}T^{13}\tau^{23}$.
\item $W:=\tau R\tau$ is a solution of the pentagonal equation.
\end{enumerate}
\end{proposition}

\begin{proof}
1 $\Leftrightarrow$ 2 The proof will follows from the
formulas:
$$
T^{12}T^{23}T^{12}=\tau^{13}R^{23}R^{13}R^{12}, \quad
T^{23}\tau^{12}T^{23}=\tau^{13}R^{12}R^{23}
$$
and from the fact that $\tau^{13}$ is an automorphism of $V\ot V\ot V$.
Let $x$, $y$, $z\in V$. Then $T(x\ot y)=\sum y_1\ot x_1$. We have
$$
T^{12}T^{23}T^{12}(x\ot y\ot z)=
\sum z_{011}\ot y_{101}\ot x_{110}=
\tau^{13}R^{23}R^{13}R^{12}(x\ot y\ot z)
$$
and
$$
T^{23}\tau^{12}T^{23}(x\ot y\ot z)=
\sum z_{10}\ot y_{11}\ot x_{01}=
\tau^{13}R^{12}R^{23}(x\ot y\ot z)
$$
1 $\Leftrightarrow$ 3 follows from the formulas:
$$
T^{23}T^{12}T^{23}=R^{23}R^{13}R^{12}\tau^{13}, \quad
T^{12}T^{13}\tau^{23}=R^{12}R^{23}\tau^{13}.
$$
1 $\Leftrightarrow$ 4 follows from the formulas:
$$
W^{12}W^{13}W^{23}=R^{23}R^{13}R^{12}\tau^{13}, \quad
W^{23}W^{12}=R^{12}R^{23}\tau^{13}.
$$
\end{proof}

From now on we shall study only the Hopf equation. In \cite{BS} numerous
examples of operators $W$ which are solutions for the pentagonal equation are
given. All these operators come from the theory of operator algebras. We
shall present only one of them, which plays a key role in classifying the
multiplicative and commutative operators (see theorem 2.2 from the above
cited paper). Let $G$ be a locally compact group and {\sl dg} a right Haar
measure on $G$. Then, $V_G(\xi)(s,t)=\xi (st,t)$ is a solution for the
pentagonal equation. It follows that $\tau V_G \tau$ is a solution for the
Hopf equation. Next, we shall present other, purely algebraic, examples of
solutions for the Hopf equation.

\begin{examples}
1. The identity map $I_{V\ot V}$ is a solution of the Hopf equation.

2. Let $V$ be a finite dimensional vector space and $u$ an automorphism
of $V$. If $R$ is a solution of the Hopf equation then
${}^{u}R:= (u\ot u)R(u\ot u)^{-1}$ is also a solution of the Hopf equation.

Indeed, as $\End_k(V\ot V)\cong \End_k(V)\ot \End_k(V)$, we can view
$R=\sum f_i\ot g_i$, where $f_i$, $g_i\in \End_k(V)$. Then
$^uR=\sum uf_iu^{-1}\ot ug_iu^{-1}$ and
$$
(^{u}R)^{12}(^{u}R)^{23}=(u\ot u\ot u)R^{12}R^{23}(u\ot u\ot u)^{-1},
$$
$$
(^{u}R)^{23}(^{u}R)^{13}(^{u}R)^{12}=
(u\ot u\ot u)R^{23}R^{13}R^{12}(u\ot u\ot u)^{-1},
$$
hence $^uR$ is also a solution of the Hopf equation.

3. Let $f$, $g\in \End_k(V)$ such that $f^2=f$, $g^2=g$
and $fg=gf$. Then, $R:=f\ot g$ is a solution of the Hopf equation.

A direct computation shows that
$$
R^{23}R^{13}R^{12}=f^2\ot fg\ot g^2, \quad  R^{12}R^{23}=f\ot gf\ot g
$$
so the above conclusion follows. With this example in mind we can
view the Hopf equation as a natural generalization of the idempotent
endomorphism. That because $R=f\ot I$ (or $R=I\ot f$) is a solution of
the Hopf equation if and only if $f^2=f$.

We suppose now that $V$ is a two dimensional vector space with
$\{v_1, v_2 \}$ a basis of $V$. Let $f_q\in \End_k(V)$ such that with
respect to the given basis is
\begin{equation}\label{proiectii}
f_q=
\left(
\begin{array}{cc}
1&q\\
0&0
\end{array}
\right)
\end{equation}
where $q$ is a scalar of $k$. Then $f_q^2=f_q$.

Let $g_q\in \End_k(V)$, $g_q=Id_V-f_q$. Then $g_q$ is also an idempotent
endomorphism of $V$ and $g_qf_q=f_qg_q$. Thus, we obtain that
$R_q=f_q\ot g_q$ with respect to the basis
$\{v_1\ot v_1, v_1\ot v_2, v_2\ot v_1, v_2\ot v_2 \}$ is
$$
R_{q}=
\left(
\begin{array}{cccc}
0&-q&0&-q^2\\
0&1&0&q\\
0&0&0&0\\
0&0&0&0
\end{array}
\right)
$$
and $R_q$ is a solution for the Hopf equation.

Now let $g=Id_V$ and $R_{q}^{\prime}=f_q\ot Id_V$. Then with respect to
the same ordonate basis of $V\ot V$, $R_{q}^{\prime}$ is given by
$$
R_{q}^{\prime}=
\left(
\begin{array}{cccc}
1&0&q&0\\
0&1&0&q\\
0&0&0&0\\
0&0&0&0
\end{array}
\right)
$$
and $R_{q}^{\prime}$ is also a solution of the Hopf equation.

4. The above $R_{q}$ and $R_{q}^{\prime}$ are also solutions of the
quantum Yang-Baxter equation, because each of them has the form $f\ot g$
with $fg=gf$. In this example we will construct a solution for the
Hopf equation which is not a solution of the quantum Yang-Baxter equation.

Let $G$ be a group and $V$ be a $G$-graded representation on $G$, that
is $V$ is a left $k[G]$-module and there exists
$\{V_{\sigma }\mid \sigma \in G \}$ a family of subspaces of $V$ such that
$$
V=\oplus_{\sigma \in G}V_{\sigma}
\quad \mbox{and} \quad g\cdot V_{\sigma}\subseteq V_{g\sigma}
$$
for all $g$, $\sigma \in G$. If $v_{\sigma}\in V_{\sigma}$ we shall write
$\mbox{deg}(v_{\sigma})=\sigma$; if $v\in V$, then $v$ is a finite sum of
homogenous elements $v=\sum v_{\sigma}$. The map
\begin{equation}\label{gr}
R:V\ot V\to V\ot V, \quad R(u\ot v)=
\sum_{\sigma}\sigma\cdot u\ot v_{\sigma}
\end{equation}
is a solution of the Hopf equation and is not a solution of the quantum
Yang-Baxter equation.

Indeed, it is enought to prove that (\ref{Heq}) holds only for homogenous
elements. Let $u_{\sigma}\in V_{\sigma}$, $u_{\tau}\in V_{\tau}$ and
$u_{\theta}\in V_{\theta}$. Then,
\begin{eqnarray*}
R^{23}R^{13}R^{12}(u_{\sigma}\ot u_{\tau}\ot u_{\theta})&=&
R^{23}R^{13}(\tau\cdot u_{\sigma}\ot u_{\tau}\ot u_{\theta})\\
&=&R^{23}(\theta\tau\cdot u_{\sigma}\ot u_{\tau}\ot u_{\theta} )\\
&=&\theta \tau \cdot u_{\sigma} \ot \theta \cdot u_{\tau} \ot u_{\theta}
\end{eqnarray*}
and
\begin{eqnarray*}
R^{12}R^{23}(u_{\sigma}\ot u_{\tau}\ot u_{\theta})&=&
R^{12}(u_{\sigma} \ot \theta \cdot u_{\tau} \ot u_{\theta})\\
\text{($\mbox{deg}(\theta\cdot u_{\tau})=\theta \tau$)}
&=&\theta \tau \cdot u_{\sigma} \ot \theta \cdot u_{\tau} \ot u_{\theta}
\end{eqnarray*}
Hence $R$ is a solution of the Hopf equation. On the other hand, by a
direct computation we get
$$
R^{12}R^{13}R^{23}(u_{\sigma}\ot u_{\tau}\ot u_{\theta})=
\theta\tau\theta\cdot u_{\sigma}\ot\theta\cdot u_{\tau}\ot u_{\theta}
$$
i.e. $R$ is not a solution of the quantum Yang-Baxter equation.

5. Let $G$ be a group and $V$ be a $G$-crossed module, that
is $V$ is a left $k[G]$-module and there exists
$\{V_{\sigma }\mid \sigma \in G \}$ a family of subspaces of $V$ such that
$$
V=\oplus_{\sigma \in G}V_{\sigma}
\quad \mbox{and} \quad g\cdot V_{\sigma}\subseteq V_{g\sigma g^{-1}}
$$
for all $g$, $\sigma \in G$. Then $R$ given by (\ref{gr}) is a solution
of the quantum Yang-Baxter equation and is not a solution of the Hopf
equation.

6. Let $q$ be a scalar of $k$, $q\neq 0$, $q\neq 1$. Then the classical
two dimensional Yang-Baxter operator
$$
R=
\left(
\begin{array}{cccc}
q&0&0&0\\
0&1&q-q^{-1}&0\\
0&0&1&0\\
0&0&0&q
\end{array}
\right)
$$
is a solution of the Yang-Baxter equation and is not a solution for the
Hopf equation.

Indeed, the element in the $(1,1)$-position of  $R^{23}R^{13}R^{12}$
is $q^3$, while the element in the $(1,1)$-position of
$R^{12}R^{23}$ is $q^2$, i.e. $R$ is not a solution of the Hopf equation.

7. Let $H$ be a bialgebra. Then
$$
R:H\ot H\to H\ot H, \quad R(g\ot h)=\sum h_{(1)}g\ot h_{(2)}
$$
for all $g$, $h\in H$, is a solution of the Hopf equation. This operator
was defined by Takesaki in \cite{Ta} for a Hopf-von Neumann algebra
$({\cal A}, \Delta)$.

8. Let $H$ be a Hopf algebra with an antipode $S$. Then $H/k$ is a
Hopf-Galois extension (see \cite{M}), i.e. the canonical map
$$
\beta:H\ot H\to H\ot H,\quad \beta(g\ot h)=\sum gh_{(1)}\ot h_{(2)}
$$
is bijective. Then $\beta$ is a solution of the Hopf equation. Furthermore,
$$
R^{\prime}: H\ot H\to H\ot H,
\quad R^{\prime}(g\ot h)=\sum g_{(1)}\ot S(g_{(2)})h
$$
is also a solution of the Hopf equation.
\end{examples}

In \cite{BS}, the concept of multiplicative and commutative (respectively
cocommutative) operator is introduced: that is, a unitary operator
$W\in {\cal L}(K\ot K)$, where $K$ is a Hilbert space, $W$ satisfies the
pentagonal equation and $W^{12}W^{23}=W^{23}W^{13}$ (respectively
$W^{12}W^{13}=W^{13}W^{12}$). We shall now introduce the corresponding
concept for the Hopf equation.

\begin{definition}
Let $V$ be a vector space and $R\in \End_k(V\ot V)$ be a solution of the
Hopf equation. Then
\begin{enumerate}
\item $R$ is called commutative if $R^{12}R^{13}=R^{13}R^{12}$.
\item $R$ is called cocommutative if $R^{13}R^{23}=R^{23}R^{13}$.
\end{enumerate}
\end{definition}

\begin{remarks}
1. Let $R\in \End_k(V\ot V)$. Then $R$ is a commutative solution of the
Hopf equation if and only if $W:=\tau R\tau$ is a commutative solution of
the pentagonal equation.

Indeed, $R^{12}R^{13}=R^{13}R^{12}$ if and only if
\begin{equation}\label{staf}
\tau^{12}W^{12}\tau^{12}\tau^{13}W^{13}\tau^{13}=
\tau^{13}W^{13}\tau^{13}\tau^{12}W^{12}\tau^{12}.
\end{equation}
Using the formulas
$$
\tau^{12}\tau^{13}=\tau^{23}\tau^{12},\quad
\tau^{13}\tau^{12}=\tau^{12}\tau^{23},
$$
$$
W^{12}\tau^{23}=\tau^{23}W^{13},\quad
\tau^{12}W^{13}=W^{23}\tau^{12},
$$
$$
W^{13}\tau^{12}=\tau^{12}W^{23},\quad
\tau^{23}W^{12}=W^{13}\tau^{23}
$$
we get that the equation (\ref{staf}) is equivalent to
$$
\tau^{12}\tau^{23}W^{13}W^{23}\tau^{12}\tau^{13}=
\tau^{13}\tau^{12}W^{23}W^{13}\tau^{23}\tau^{12}.
$$
The conclusion follows as
$$
\tau^{12}\tau^{13}\tau^{12}\tau^{23}=
\tau^{23}\tau^{12}\tau^{13}\tau^{12}=Id.
$$
2. Suppose that $R\in \End_k(V\ot V)$ is bijective. Then, $R$ is a
cocomutative solution of the Hopf equation if and only if
$\tau R^{-1}\tau$ is a commutative solution of the Hopf equation.
\end{remarks}

Our example (4) can be generalized to arbitrary Hopf modules and
evidences the role which can be played by the $H$-Hopf modules in solving
the Hopf equation.

\begin{proposition}
Let $H$ be a bialgebra and $(M,\cdot, \rho)$ an $H$-Hopf module.
Then:
\begin{enumerate}
\item the natural map
$$R_{(M,\cdot ,\rho )}(m\ot n)=\sum n_{<1>}\cdot m\ot n_{<0>}$$
is a solution of the Hopf equation.
\item if $H$ is commutative then $R_{(M,\cdot ,\rho )}$ is a commutative
solution of the Hopf equation.
\end{enumerate}
\end{proposition}

\begin{proof}
1. Let $R=R_{(M,\cdot ,\rho )}$. For $l$, $m$, $n\in M$ we have
\begin{eqnarray*}
R^{12}R^{23}(l\ot m\ot n)&=&
R^{12}\Bigl ( \sum l\ot n_{<1>}\cdot m \ot n_{<0>} \Bigl )\\
&=&\sum (n_{<1>}\cdot m)_{<1>}\cdot l \ot (n_{<1>}\cdot m)_{<0>}\ot n_{<0>}
\end{eqnarray*}
and
\begin{eqnarray*}
R^{23}R^{13}R^{12}(l\ot m\ot n)&=&
R^{23}R^{13}\Bigl (\sum m_{<1>}\cdot l\ot m_{<0>}\ot n \Bigl)\\
&=&R^{23}\Bigl ( \sum n_{<1>}m_{<1>}\cdot l \ot m_{<0>}\ot n_{<0>}\Bigl )\\
&=&\sum n_{<2>}m_{<1>}\cdot l\ot n_{<1>}\cdot m_{<0>}\ot n_{<0>}\\
&=&\sum n_{<1>(2)}m_{<1>}\cdot l\ot n_{<1>(1)}\cdot m_{<0>}\ot n_{<0>}\\
\text{(\mbox{using} (\ref{H1}))}
&=&\sum (n_{<1>}\cdot m)_{<1>}\cdot l \ot (n_{<1>}\cdot m)_{<0>}\ot n_{<0>}
\end{eqnarray*}
i.e. $R$ is a solution of the Hopf equation.

2. We have
$$
R^{12}R^{13}(l\ot m\ot n)=\sum m_{<1>}n_{<1>}\cdot l \ot m_{<0>} \ot n_{<0>}
$$
and
$$
R^{13}R^{12}(l\ot m\ot n)=\sum n_{<1>}m_{<1>}\cdot l \ot m_{<0>} \ot n_{<0>}
$$
As $H$ is commutative, we obtain that $R^{12}R^{13}=R^{13}R^{12}$.
\end{proof}

\begin{remark}
If $(M,\cdot, \rho)$ is an $H$-Hopf module then the map
$$
R_{(M,\cdot ,\rho )}^{\prime}:M\ot M\to M\ot M, \quad
R_{(M,\cdot ,\rho )}^{\prime}(m\ot n)=\sum m_{<0>}\ot m_{<1>}\cdot n
$$
is a solution of the pentagonal equation, as
$R_{(M,\cdot ,\rho )}^{\prime}=\tau R_{(M,\cdot ,\rho )}\tau $.
\end{remark}

\section{A FRT type construction for Hopf modules}

In this section we shall prove the main result of the paper, which
shows us that in the finite dimensional case any solution of the Hopf
equation  has the form $R_{(M,\cdot, \rho)}$.

\begin{theorem}
Let $M$ be a finite dimensional vector space and $R\in \End_k(M\ot M)$
be a solution of the Hopf equation. Then
\begin{enumerate}
\item There exists a bialgebra $B(R)$ such that $M$ has a structure
of $B(R)$-Hopf module $(M,\cdot, \rho)$ and $R=R_{(M,\cdot, \rho)}$.
\item The bialgebra $B(R)$ is a universal object with this property:
if $H$ is a bialgebra such that
$(M,\cdot^{\prime}, \rho^{\prime})\in {}_H{\cal M}^H$ and
$R=R_{(M,\cdot^{\prime}, \rho^{\prime})}$ then there exists a unique
bialgebra map $f:B(R)\to H$ such that $\rho^{\prime}=(I\ot f)\rho$.
Furthermore, $a\cdot m=f(a)\cdot^{\prime}m$, for all $a\in B(R)$,
$m\in M$.
\item If $R$ is commutative, then there exists a commutative bialgebra
$\overline{B}(R)$ such that $M$ has a structure
of $\overline{B}(R)$-Hopf module $(M,\cdot^{\prime}, \rho^{\prime})$ and
$R=R_{(M,\cdot^{\prime}, \rho^{\prime})}$.
\end{enumerate}
\end{theorem}

\begin{proof}
$1$. The proof will be given is several steps. Let $\{m_1,\cdots ,m_n \}$
be a basis for $M$ and $(x_{uv}^{ji})_{i,j,u,v}$ a family
of scalars of $k$ such that
\begin{equation}
R(m_v\ot m_u)=\sum_{i,j}x_{uv}^{ji}m_i\ot m_j
\end{equation}
for all $u$, $v=1,\cdots ,n$.

Let $(C, \Delta, \varepsilon)={\cal M}^n(k)$, be the comatrix coalgebra
of order $n$, i.e. $C$ is the coalgebra with the basis
$\{c_{ij}\mid i,j=1,\cdots,n \}$ such that
\begin{equation}\label{com1}
\Delta(c_{jk})=\sum_{u=1}^{n}c_{ju}\ot c_{uk},
\quad \varepsilon (c_{jk})=\delta_{jk}
\end{equation}
for all $j,k=1,\cdots, n$. Let $\rho :M\to M\ot C$ given by
\begin{equation}\label{ro}
\rho(m_l)=\sum_{v=1}^{n}m_v\ot c_{vl}
\end{equation}
for all $l=1,\cdots, n$. Then, by lemma \ref{unu}, $M$ is a right
$C$-comodule. Let $T(C)$ be the bialgebra structure on the tensor algebra
$T(C)$ which extends $\Delta$ and $\varepsilon$ (from lemma \ref{doi}).
As the inclusion $i:C\to T(C)$ is a coalgebra map, $M$ has a right
$T(C)$-comodule structure via
$$
M\stackrel{\rho}{\longrightarrow}
M\ot C\stackrel{I\ot i}{\longrightarrow}M\ot T(C)
$$
There will be no confusion if we also denote the right
$T(C)$-comodule structure on $M$ with $\rho$.

Now, we will put a left $T(C)$-module structure on $M$ in such a way that
$R=R_{(M,\cdot, \rho)}$. First we define
$$
\mu :C\ot M\to M, \quad
\mu(c_{ju}\ot m_v):=\sum_{i} x_{uv}^{ji}m_i
$$
for all $j$, $u$, $v=1,\cdots n$. From lemma \ref{doi}, there exists
a unique left $T(C)$-module structure on $(M,\cdot )$ such that
$$
c_{ju}\cdot m_v=\mu(c_{ju}\ot m_v)=\sum_{i} x_{uv}^{ji}m_i
$$
for all $j$, $u$, $v=1,\cdots ,n$. For $m_v$, $m_u$ the elements of
the given basis, we have:
\begin{eqnarray*}
R_{(M,\cdot, \rho)}(m_v\ot m_u)&=&
\sum_j c_{ju}\cdot m_v\ot m_j\\
&=&\sum_{i,j} x_{uv}^{ji}m_i\ot m_j\\
&=&R(m_v\ot m_u)
\end{eqnarray*}
Hence, $(M,\cdot, \rho)$ has a structure of left $T(C)$-module and
right $T(C)$-comodule such that $R=R_{(M,\cdot, \rho)}$.

Now, we define the {\sl obstructions} $\chi (i,j,k,l)$ which measure
how far away $M$ is from a $T(C)$-Hopf module. Keeping in mind that
$T(C)$ is generated as an algebra by $(c_{ij})$ and using lemma
\ref{trei} we compute
$$
\sum h_{(1)}\cdot m_{<0>}\ot h_{(2)}m_{<1>} - \rho (h\cdot m)
$$
only for $h=c_{jk}$, and $m=m_l$, for $j$, $k$, $l=1,\cdots, n$. We have:
\begin{eqnarray*}
\sum h_{(1)}\cdot m_{<0>}\ot h_{(2)}m_{<1>}&=&
\sum_{u,v}c_{ju}\cdot m_v\ot c_{uk}c_{vl}\\
&=&\sum_{i}m_i \ot \Bigl(\sum_{u,v}x_{uv}^{ji}c_{uk}c_{vl} \Bigl)
\end{eqnarray*}
and
\begin{eqnarray*}
\rho(h\cdot m)&=&\rho (c_{jk}\cdot m_l)\\
&=&\sum_{\alpha}x_{kl}^{j\alpha}(m_{\alpha})_{<0>}\ot (m_{\alpha})_{<1>}\\
&=&\sum_{i,\alpha}x_{kl}^{j\alpha}m_i\ot c_{i\alpha}\\
&=&\sum_im_i\ot \Bigl(\sum_{\alpha}x_{kl}^{j\alpha}c_{i\alpha}\Bigl)
\end{eqnarray*}
Let
\begin{equation}\label{obst}
\chi(i,j,k,l):=\sum_{u,v}x_{uv}^{ji}c_{uk}c_{vl} -
\sum_{\alpha}x_{kl}^{j\alpha}c_{i\alpha}
\end{equation}
for all $i$, $j$, $k$, $l=1,\cdots, n$. Then
\begin{equation}\label{pais}
\sum h_{(1)}\cdot m_{<0>}\ot h_{(2)}m_{<1>}-
\rho(h\cdot m)=\sum_im_i\ot \chi(i,j,k,l)
\end{equation}
Let $I$ be the two-sided ideal of $T(C)$ generated by all $\chi(i,j,k,l)$,
$i$, $j$, $k$, $l=1,\cdots, n$. The key point of
the proof follows:

{\sl $I$ is a bi-ideal of $T(C)$ and $I\cdot M=0$.}

We first prove that $I$ is also a coideal and this will result from the
following formula:
\begin{equation}\label{prima}
\Delta(\chi(i,j,k,l))=\sum_{a,b}\chi(i,j,a,b)\ot c_{ak}c_{bl}+
\sum_{p}c_{ip}\ot \chi(p,j,k,l)
\end{equation}
Indeed, we have:
\begin{eqnarray*}
\Delta(\chi(i,j,k,l))&=&\sum_{u,v}x_{uv}^{ji}\Delta(c_{uk})\Delta(c_{vl})-
\sum_{\alpha}x_{kl}^{j\alpha}\Delta(c_{i\alpha})\\
&=&\sum_{a,b,u,v}x_{uv}^{ji}c_{ua}c_{vb}\ot c_{ak}c_{bl}-
\sum_{\alpha, p}x_{kl}^{j\alpha}c_{ip}\ot c_{p\alpha}\\
&=&\sum_{a,b}\Bigl(\sum_{u,v}x_{uv}^{ji}c_{ua}c_{vb}\Bigl)
\ot c_{ak}c_{bl} -
\sum_{p}c_{ip}\ot \Bigl( \sum_{\alpha}x_{kl}^{j\alpha}c_{p\alpha}\Bigl)\\
&=&\sum_{a,b}\Bigl(\chi(i,j,a,b)+
\sum_{\gamma}x_{ab}^{j\gamma}c_{i\gamma}\Bigl)\ot c_{ak}c_{bl}\\
&-&\sum_{p}c_{ip}\ot \Bigl(-\chi(p,j,k,l)+
\sum_{r,s}x_{rs}^{jp}c_{rk}c_{sl}\Bigl) \\
&=&\sum_{a,b}\chi(i,j,a,b)\ot c_{ak}c_{bl}+
\sum_{p}c_{ip}\ot \chi(p,j,k,l)
\end{eqnarray*}
where in the last equality we use the fact that
$$
\sum_{a,b,\gamma}x_{ab}^{j\gamma}c_{i\gamma}\ot c_{ak}c_{bl}=
\sum_{p,r,s}x_{rs}^{jp}c_{ip}\ot c_{rk}c_{sl}
$$
Hence, the formula (\ref{prima}) holds. On the other hand
$$
\varepsilon\Bigl(\chi(i,j,k,l)\Bigl)=x_{kl}^{ji}-x_{kl}^{ji}=0
$$
so we proved that $I$ is a coideal of $T(C)$.

Now, in order to show that $I\cdot M=0$, we shall use the fact that
$R$ is a solution of the Hopf equation. For $z\in M$, $j$, $k=1,\cdots ,n$,
we have the following formula:
\begin{equation}\label{adoua}
\Bigl(R^{23}R^{13}R^{12}-R^{12}R^{23}\Bigl)(z\ot m_k\ot m_j)=
\sum_{r,s}\chi(r,s,j,k)\cdot z\ot m_r\ot m_s
\end{equation}
Let us compute
\begin{eqnarray*}
\Bigl(R^{23}R^{13}R^{12} \Bigl)(z\ot m_k\ot m_j)&=&
\Bigl(R^{23}R^{13}\Bigl)
(\sum_{\alpha}c_{\alpha k}\cdot z\ot m_{\alpha}\ot m_{j} )\\
&=&R^{23}\Bigl(\sum_{\alpha,\beta}
c_{\beta j}c_{\alpha k}\cdot z\ot m_{\alpha}\ot m_{\beta} ) \Bigl)\\
&=&\sum_{\alpha, \beta, r,s}x_{\beta\alpha}^{sr}
c_{\beta j}c_{\alpha k}\cdot z\ot m_{r}\ot m_{s}
\end{eqnarray*}
On the other hand
\begin{eqnarray*}
\Bigl(R^{12}R^{23}\Bigl)(z\ot m_k\ot m_j)&=&
R^{12}(\sum_{s}z\ot c_{sj}\cdot m_k \ot m_s)\\
&=&R^{12}(\sum_{s,\alpha}z\ot x_{jk}^{s\alpha}m_{\alpha}\ot m_s)\\
&=&\sum_{r,s,\alpha}x_{jk}^{s\alpha}c_{r\alpha}\cdot z\ot m_r\ot m_s
\end{eqnarray*}
It follows that
\begin{eqnarray*}
\Bigl(R^{23}R^{13}R^{12}-R^{12}R^{23}\Bigl)(z\ot m_k\ot m_j)
&=&\sum_{r,s}\Bigl(\sum_{\alpha,\beta}
x_{\beta\alpha}^{sr}c_{\beta j}c_{\alpha k}-
\sum_{\alpha}x_{jk}^{s\alpha}c_{r\alpha} \Bigl)\cdot z\ot m_r\ot m_s\\
&=&\sum_{r,s}\chi(r,s,j,k)\cdot z\ot m_r\ot m_s
\end{eqnarray*}
i.e. the formula (\ref{adoua}) holds. But $R$ is a solution of the
Hopf equation, hence $\chi(r,s,j,k)\cdot z =0$, for all $z\in M$,
$j$, $k$, $r$, $s=1,\cdots, n$. We conclude that
$I$ is a bi-ideal of $T(C)$ and $I\cdot M=0$. Define now
$$
B(R)=T(C)/I.
$$
$M$ has a right $B(R)$-comodule structure via the canonical projection
$T(C)\to B(R)$ and a left $B(R)$-module structure as $I\cdot M=0$.
As $(c_{ij})$ generate $B(R)$ and in $B(R)$, $\chi(i,j,k,l)=0$, for all
$i$, $j$, $k$, $l=1,\cdots, n$, using (\ref{pais}) we get that
$(M,\cdot, \rho)\in {}_{B(R)}{\cal M}^{B(R)}$ and, by lemma (\ref{patru}),
$R=R_{(M,\cdot,\rho)}$.

$2$. Let $H$ be a bialgebra and suppose that
$(M,\cdot^{\prime}, \rho^{\prime})\in {}_H{\cal M}^H$ and
$R=R_{(M,\cdot^{\prime}, \rho^{\prime})}$.
Let $(c_{ij}^{\prime})_{i,j=1,\cdots ,n}$ be a family of elements of $H$
such that
$$
\rho^{\prime}(m_l)=\sum_v m_v\ot c_{vl}^{\prime}
$$
Then
$$
R(m_v\ot m_u)=\sum_j c_{ju}^{\prime}\cdot^{\prime}m_v \ot m_j
$$
and
$$
c_{ju}^{\prime}\cdot^{\prime}m_v=\sum_i x_{uv}^{ji}m_i=c_{ju}\cdot m_v.
$$
Let
$$
\chi^{\prime}(i,j,k,l)=
\sum_{u,v}x_{uv}^{ji}c_{uk}^{\prime}c_{vl}^{\prime} -
\sum_{\alpha}x_{kl}^{j\alpha}c_{i\alpha}^{\prime}
$$
From the universal property of the tensor algebra $T(C)$, there exists
a unique algebra map $f_1: T(C)\to H$ such that
$f_1(c_{ij})=c_{ij}^{\prime}$, for all $i$, $j=1,\cdots, n$.
As $(M,\cdot^{\prime}, \rho^{\prime})\in {}_H{\cal M}^H$ we get that
$\chi^{\prime}(i,j,k,l)=0$, and hence $f_1(\chi(i,j,k,l))=0$,
for all $i$, $j$, $k$, $l=1,\cdots, n$. So the map $f_1$ factorizes to
the map
$$
f:B(R)\to H, \quad f(c_{ij})=c_{ij}^{\prime}
$$
Of course, for $m_l$ an arbitrary element of the given basis of $M$, we have
$$
(I\ot f)\rho(m_l)=\sum_v m_v\ot f(c_{vl})=\sum_v m_v\ot c_{vl}^{\prime}=
\rho^{\prime}(m_l)
$$
Conversely, the relation $(I\ot f)\rho=\rho^{\prime}$ necessarily
implies $f(c_{ij})=c_{ij}^{\prime}$, which proves the uniqueness of $f$.
This completes the proof of the theorem.

3. For $z\in M$ and $j,k=1,\cdots, n$ we have the formula:
\begin{equation}\label{dorm}
\Bigl(R^{12}R^{13}-R^{13}R^{12}\Bigl)(z\ot m_k\ot m_j)=
\sum_{r,s}\Bigl(c_{rk}c_{sj}-c_{sj}c_{rk}\Bigl)\cdot z\ot m_r\ot m_s
\end{equation}
Let $\overline{I}$ be the two-sided ideal of $T(C)$ generated by $I$ and all
$[c_{rk}, c_{sj}]$. It follows from the formula
$$
\Delta\Bigl([c_{rk}, c_{sj}]\Bigl)=\sum_{a,b} \Bigl([c_{ra}, c_{sb}]\ot
c_{bj}c_{ak}+ c_{ra}c_{sb}\ot [c_{ak}, c_{bj}] \Bigl)
$$
that $\overline{I}$ is also a coideal of $T(C)$ and from equation
(\ref{dorm}) we get that $\overline{I}\cdot M=0$. Define now
$$
\overline{B}(R)=T(C)/ \overline{I}.
$$
Then $\overline{B}(R)$ is a commutative bialgebra, $M$ has a structure
of $\overline{B}(R)$-Hopf module $(M,\cdot^{\prime}, \rho^{\prime})$ and
$R=R_{(M,\cdot^{\prime}, \rho^{\prime})}$.
\end{proof}

\begin{remark}
1. Our obstruction elements $\chi(i,j,k,l)$ play the same role as the
homogenous elements $d(i,j,k,l)$ defined in \cite{R1} which correspond
to the quantum Yang-Baxter equation: the two-sided ideal generated by
them is also a coideal which anihilates $M$. This was the key point of
the proof.

2. The last point of our theorem can by viewed as an algebraic version of
theorem 2.2. from \cite{BS}. All commutative bialgebras $\overline{B}(R)$
are quotients for various bialgebra structures which can be given on
$k[Y_1,\cdots, Y_n]$.
\end{remark}

\section{Applications}
In this section we shall construct new examples of noncommutative
noncocommutative bialgebras arising from our FRT type theorem.
As the relations through which we factor are not all homogenous, all our
examples are different from the ones which appear in quantum group
theory. A completely different method for constructing such objects
uses Ore extensions and was recently evidenced in \cite{BDG}.

In the next propositions of this section, the relations $\chi(i,j,k,l)=0$
will be written in the lexicografical order according to $(i,j,k,l)$
starting with $(1,1,1,1)$.

\subsection{Back to euclidian geometry}

There exists an intimate link between the quantum Yang-Baxter equation and
the quantum plane $k_{q}<x, y \mid xy=qyx>$ (see \cite{K}).

We shall now show that the term "dequantisation" used in the introduction is
not an abuse. More specifically, our FRT type construction for the Hopf
equation supplies us with a way to construct noncommutative and
noncocommutative bialgebras starting from the elementary maps of plane
euclidian geometry: projections of $k^2$ on the $Ox$ and $Oy$ axis.
The map $f_0$, corresponding to $q=0$ in the equation (\ref{proiectii}),
is in fact the projection of the plane $k^2$ on the $Ox$ axis.
We can now associate three bialgebras to
this projection: the first one corresponds to the solution of the Hopf
equation $f_0\ot g_0$, where $g_0=Id_{k^{2}}-f_0$,
i.e. is the projection  of the plane $k^2$ on the $Oy$ axis,
the second corresponds to $f_0\ot Id_{k^2}$, and the third correspunds
to the $f_0\ot f_0$.
In this way, we obtain the bialgebras denoted below by $B_{0}^{2}(k)$,
$D_{0}^{2}(k)$ and $E_{0}^{2}(k)$.
If $q\neq 0$ and $k=\mathbf{R}$, the map $f_q$ given in (\ref{proiectii})
sends all the points of the $\mathbf{R}^2$ plane on the $Ox$ axis
under an angle arctg($q$) with respect to the $Oy$ axis.
Correspondingly, the bialgebras $B_{q}^{2}(k)$, $D_{q}^{2}(k)$ and
$E_{q}^{2}(k)$ are constructed.

\begin{proposition} \label{ox}
Let $q$ be a scalar of the field $k$ and $R_q$ be the solution of the
Hopf equation given by
$$
R_{q}=
\left(
\begin{array}{cccc}
0&-q&0&-q^2\\
0&1&0&q\\
0&0&0&0\\
0&0&0&0
\end{array}
\right)
$$
Let $B_q^{2}(k)$ be the bialgebra $B(R_q)$. Then:
\begin{enumerate}
\item If $q=0$, the bialgebra $B_0^{2}(k)$ is the free algebra generated
by $x$, $y$, $z$ with the relations:
$$
yx=x, \quad yz=0.
$$
The comultiplication $\Delta$ and the counity $\varepsilon$ are given by:
$$
\Delta(x)=x\ot x,\quad \Delta(y)=y\ot y,\quad \Delta(z)=x\ot z+z\ot y
$$
$$
\varepsilon(x)=\varepsilon(y)=1,\quad \varepsilon(z)=0.
$$
\item If $q\neq 0$, the bialgebra $B_q^{2}(k)$ is the free algebra generated
by $A$, $B$ with the relation:
$$
A^2B=AB.
$$
The comultiplication $\Delta$ and the counity $\varepsilon$ are given by:
$$
\Delta(A)=A\ot A,\quad \Delta(B)=q^{-1}AB\ot B+(B-AB)\ot A,
$$
$$
\varepsilon(A)=1,\quad \varepsilon(B)=q.
$$
\end{enumerate}
\end{proposition}
\begin{proof}
Let $M$ be a two dimensional vector space with $\{m_1,m_2 \}$ a basis.
Put $R=R_q$. With respect to the ordonate basis
$\{m_1\ot m_1, m_1\ot m_2, m_2\ot m_1, m_2\ot m_2 \}$, $R$ is given by:
$$
R(m_1\ot m_1)=R(m_2\ot m_1)=0,
$$
$$
R(m_1\ot m_2)=-qm_1\ot m_1 +m_1\ot m_2,\quad
R(m_2\ot m_2)=-q^2m_1\ot m_1 +qm_1\ot m_2
$$
Now, if we write
$$
R(m_v\ot m_u)=\sum_{i,j=1}^{2}x_{uv}^{ji}m_i\ot m_j
$$
we get that among the elements $(x_{uv}^{ji})$, the only nonzero elements
are:
$$
x_{21}^{11}=-q,\quad x_{21}^{21}=1,\quad
x_{22}^{11}=-q^2,\quad x_{22}^{21}=q.
$$
The sixteen relations $\chi(i,j,k,l)=0$ are:
$$
-qc_{21}c_{11}-q^2c_{21}c_{21}=0,\quad
-qc_{21}c_{12}-q^2c_{21}c_{22}=0,
$$
$$
-qc_{22}c_{11}-q^2c_{22}c_{21}=-qc_{11},\quad
-qc_{22}c_{12}-q^2c_{22}c_{22}=-q^2c_{11},
$$
$$
c_{21}c_{11}+qc_{21}c_{21}=0,\quad
c_{21}c_{12}+qc_{21}c_{22}=0,
$$
$$
c_{22}c_{11}+qc_{22}c_{21}=c_{11}, \quad
c_{22}c_{12}+qc_{22}c_{22}=qc_{11},
$$
$$
0=0, \quad 0=0, \quad 0=-qc_{21}, \quad 0=-q^2c_{21}
$$
$$
0=0, \quad 0=0, \quad 0=c_{21}, \quad 0=qc_{21}.
$$
Hence, $c_{21}=0$. Now, if we denote $c_{11}=x$, $c_{22}=y$, $c_{12}=z$,
there are only two linear independent relations:
$$
yx=x, \quad yz+qy^2=qx.
$$
If $q=0$, then follows 1. If $q\neq 0$, then $x$ is an element in the
free algebra generated by $y$ and $z$. Let $A=y$ and $B=z+qy=z+qA$. Then
$$
x=q^{-1}AB
$$
and by substituting in the first relation we get $A^2B=AB$. The formulas for
$\Delta$ and $\varepsilon$ follow as the original $(c_{ij})$ was
a comultiplicative matrix.
\end{proof}

\begin{remark}
The bialgebra $B_0^{2}(k)$ is not a Hopf algebra. We can localize it to
obtain a Hopf algebra. As $\Delta(x)=x\ot x$, $\Delta(y)=y\ot y$ and
$\varepsilon(x)=\varepsilon(y)=1$ we should add new generators which
make $x$ and $y$ invertible. But then $y=1$ and $z=0$. It
follows that if we localize the bialgebra $B_0^{2}(k)$, we get the usual
Hopf algebra $k[X,X^{-1}]$, with $\Delta (X)=X\ot X$,
$\varepsilon(X)=X\ot X$, and with the antipode $S(X)=X^{-1}$.
\end{remark}

\begin{proposition}
Let $q$ be a scalar of the field $k$ and $R_q^{\prime}$ be the
solution of the Hopf equation given by
$$
R_{q}^{\prime}=
\left(
\begin{array}{cccc}
1&0&q&0\\
0&1&0&q\\
0&0&0&0\\
0&0&0&0
\end{array}
\right)
$$
We denote by $D_{q}^{2}(k)$ the bialgebra $B(R_q^{\prime})$. Then:
\begin{enumerate}
\item If $q=0$, the bialgebra $D_0^{2}(k)$ is the free algebra generated by
$x$, $y$, $z$ with the relations:
$$
x^2=x=yx, \quad zx=xz=z^2=yz=0.
$$
The comultiplication $\Delta$ and the counity $\varepsilon$ are given by:
$$
\Delta(x)=x\ot x,\quad \Delta(y)=y\ot y,\quad \Delta(z)=x\ot z+z\ot y
$$
$$
\varepsilon(x)=\varepsilon(y)=1,\quad \varepsilon(z)=0.
$$
\item If $q\neq 0$, the bialgebra $D_q^{2}(k)$ is the free algebra generated
by $A$, $B$ with the relations:
$$
A^3=A^2,\qquad BA=0.
$$
The comultiplication $\Delta$ and the counity $\varepsilon$ are given by:
$$
\Delta(A)=A\ot A +q^{-1}(A^2-A)\ot B, \quad
\Delta(B)=A^2\ot B +B\ot A -q^{-1}B\ot B,
$$
$$
\varepsilon(A)=1, \quad \varepsilon(B)=0.
$$
\end{enumerate}
\end{proposition}

\begin{proof}
We start exactly as in the above proposition. We get that among the
scalars $(x_{uv}^{ji})$, which define $R$, the only nonzero elements are:
$$
x_{11}^{11}=x_{21}^{21}=1, \quad x_{12}^{11}=x_{22}^{21}=q.
$$
Now the relations $\chi(i,j,k,l)=0$ are:
$$
c_{11}c_{11}+qc_{11}c_{21}=c_{11}, \quad
c_{11}c_{12}+qc_{11}c_{22}=qc_{11},
$$
$$
c_{12}c_{11}+qc_{12}c_{21}=0, \quad
c_{12}c_{12}+qc_{12}c_{22}=0,
$$
$$
c_{21}c_{11}+qc_{21}c_{21}=0, \quad
c_{21}c_{12}+qc_{21}c_{22}=0,
$$
$$
c_{22}c_{11}+qc_{22}c_{21}=c_{11}, \quad
c_{22}c_{12}+qc_{22}c_{22}=qc_{11},
$$
$$
0=c_{21}, \quad 0=qc_{21}, \quad 0=0, \quad 0=0,
$$
$$
0=0, \quad 0=0, \quad 0=c_{21}, \quad 0=qc_{21}.
$$
Hence $c_{21}=0$. If we denote $c_{11}=x$, $c_{22}=y$, $c_{12}=z$ then
we get the following six relations:
$$
x^2=x=yx, \quad zx=0, \quad z^2+qzy=0,
$$
$$
xz+qxy=yz+qy^2=qx.
$$
As $c_{21}=0$, the comultiplication and the counity take the
following form:
\begin{equation}\label{dudu}
\Delta(x)=x\ot x,\quad \Delta(y)=y\ot y,\quad \Delta(z)=x\ot z+z\ot y
\end{equation}
$$
\varepsilon(x)=\varepsilon(y)=1,\quad \varepsilon(z)=0.
$$
Hence, if $q=0$ we get exactly the relation of $D_{0}^{2}(k)$. Suppose now
that $q\neq 0$. Then $x$ is in the free algebra generated by $y$, $z$
and
$$
x=y^2 +q^{-1}yz=y(y +q^{-1}z)
$$
If we substitute $x$ in the other five relations we get, after we
multiply with $q$ or $q^2$:
$$
y(z+qy)[y(z+qy)-q]=0,
$$
$$
(y-1)y(z+qx)=0,
$$
$$
zy(z+qy)=0,
$$
$$
z(z+qy)=0,
$$
$$
y(z+qy)(z+qy-q)=0.
$$
Using the fact that $y(y-1)=(y-1)y$, the fifth relation follows from the
second and the fourth. On the other hand, as $y(y^2-1)=(y+1)(y-1)y$, the
first relation follows from the second and the third. We have thus reduced
the above five relations to only three:
$$
(y-1)y(z+qy)=0,\quad zy(z+qy)=0,\quad z(z+qy)=0.
$$
Further, let $B=z$ and $A=q^{-1}(z+qy)=q^{-1}(B+qy)$.
It follows that $y=A-q^{-1}B$. The third relation takes the form
$$BA=0$$
and this implies the second one.
Using $BA=0$, the first relation becomes $A^3=A^2$.
The formulas for $\Delta$ and
$\varepsilon$ follow from equation (\ref{dudu}). We note that the element
$A-q^{-1}B$ is a groupal element of $D_{q}^{2}(k)$.
\end{proof}

\begin{remark}
$R_{q}^{\prime}$ is also a solution of the quantum Yang-Baxter equation.
The bialgebra $A(R_{q}^{\prime})$ which we obtain applying the usual FRT
construction is the free algebra generated by $x$, $y$, $z$, $t$ with the
relations:
$$
zx=xz=zy=z^2=zt=0, \quad xy-yx=qyz, \quad xt-tx=qtz,
$$
$$
xy+qxt=qx^2, \quad y^2+qyt=qxy, \quad ty+qt^2=qxt.
$$
The comultiplication $\Delta$ and the counity $\varepsilon$ are given
in such way that the matrix
$$
\left(
\begin{array}{cc}
x&y\\
z&t
\end{array}
\right)
$$
is comultiplicative.
\end{remark}

In the next proposition we shall prove that the bialgebra $E_{q}^{2}(k)$,
with $q\neq 0$, which can be associated to the solution
$R_{q}^{\prime\prime}=f_q\ot f_q$ is not dependent of $q$, i.e.
$E_{q}^{2}(k)\cong E_{q^{\prime}}^{2}(k)$, for all $q$,
$q^{\prime}\in k\backslash \{0\}$.

\begin{proposition}
Let $q$ be a scalar of the field $k$ and $R_q^{\prime\prime}$ be the
solution of the Hopf equation given by
$$
R_{q}^{\prime\prime}=
\left(
\begin{array}{cccc}
1&q&q&q^2\\
0&0&0&0\\
0&0&0&0\\
0&0&0&0
\end{array}
\right)
$$
Let $E_q^{2}(k)$ be the bialgebra $B(R_{q}^{\prime\prime})$. Then:
\begin{enumerate}
\item If $q=0$, the bialgebra $E_0^{2}(k)$ is the free algebra generated
by $x$, $y$, $z$ with the relations:
$$
x^2=x, \quad xz=zx=z^2=0.
$$
The comultiplication $\Delta$ and the counity $\varepsilon$ are given by:
$$
\Delta(x)=x\ot x,\quad \Delta(y)=y\ot y,\quad \Delta(z)=x\ot z+z\ot y
$$
$$
\varepsilon(x)=\varepsilon(y)=1,\quad \varepsilon(z)=0.
$$
\item If $q\neq 0$, the bialgebra $E_q^{2}(k)$ is the free algebra generated
by $A$, $B$ with the relations:
$$
B^{3}=B^{2}.
$$
The comultiplication $\Delta$ and the counity $\varepsilon$ are given by:
$$
\Delta(A)=A\ot A,\quad \Delta(B)=B\ot A +B^{2}\ot (B-A)
$$
$$
\varepsilon(A)=\varepsilon(B)=1.
$$
\end{enumerate}
\end{proposition}

\begin{proof}
We start exactly as in the above propositions. We get that among the
scalars $(x_{uv}^{ji})$, which define $R$, the only nonzero elements are:
$$
x_{11}^{11}=1, \quad x_{21}^{11}=x_{12}^{11}=q, \quad x_{22}^{11}=q^2.
$$
Now the relations $\chi(i,j,k,l)=0$ are:
$$
c_{11}c_{11}+qc_{21}c_{11}+qc_{11}c_{21}+q^2c_{21}c_{21}=c_{11}
$$
$$
c_{11}c_{12}+qc_{21}c_{12}+qc_{11}c_{22}+q^2c_{21}c_{22}=qc_{11}
$$
$$
c_{12}c_{11}+qc_{22}c_{11}+qc_{12}c_{21}+q^2c_{22}c_{21}=qc_{11}
$$
$$
c_{12}c_{12}+qc_{22}c_{12}+qc_{12}c_{22}+q^2c_{22}c_{22}=q^2c_{11}
$$
$$
0=0, \quad 0=0, \quad 0=0,\quad 0=0,
$$
$$
0=c_{21}, \quad 0=qc_{21}, \quad 0=qc_{21},\quad 0=q^2c_{21},
$$
$$
0=0, \quad 0=0, \quad 0=0,\quad 0=0,
$$
Hence $c_{21}=0$. If we denote $c_{11}=x$, $c_{22}=y$, $c_{12}=z$ then
we get the following four relations:
$$
x^2=x,\quad xz+qxy=qx,\quad zx+qyx=qx,
$$
$$
z^2+qyz+qzy+q^2y^2=q^2x.
$$
So, if $q=0$, we obtain the relations of $E_{0}^{2}(k)$. If $q\neq 0$, then
$x$ is in the free algebra generated by $y$ and $z$ and
$$
x=y^2+q^{-1}zy+q^{-1}yz+q^{-2}z^2=(y+q^{-1}z)^{2}.
$$
If we substitute $x$ in the other three relations, only
$$
(y+q^{-1}z)^{3}=(y+q^{-1}z)^{2}.
$$
remains, the other two being linear dependent from this one.
Now, if we denote $A=y$ and $B=y+q^{-1}z$, we obtain the description of
$E_{q}^{2}(k)$.
\end{proof}

\begin{remarks}
1. If $q\neq 0$, it is interesting to denote that $B^2$ is a groupal element
of $E_{q}^{2}(k)$. Indeed, we have
\begin{eqnarray*}
\Delta(B^2)&=&B^2\ot A^2+B^2\ot (AB-A^2)+B^2\ot (BA-A^2)+
B^4\ot (B^2-BA-AB+A^2)\\
&=&B^2\ot A^2+B^2\ot (AB-A^2)+B^2\ot (BA-A^2)+
B^2\ot (B^2-BA-AB+A^2)\\
&=&B^2\ot B^2
\end{eqnarray*}
2. The bialgebra $D_{0}^{2}(k)$ is the quotient of $B_{0}^{2}(k)$
by the two-sided ideal (which is also a coideal) generated by
$$
x^2-x, \quad zx, \quad xz, \quad z^2.
$$
$D_{0}^{2}(k)$ is also a quotient of $E_{0}^{2}(k)$ by the two-sided
ideal generated by
$$
yx-x, \quad yz.
$$
3. Let $n\geq 2$ be a natural number. The bialgebras $B_{0}^{2}(k)$,
$D_{0}^{2}(k)$ and $E_{0}^{2}(k)$ constructed in the previous propositions
can be generalised to
$B_{0}^{n}(k)$, $D_{0}^{n}(k)$ and respectively $E_{0}^{n}(k)$.
We have chosen to construct them for the case $n=2$ in order to better sense
the flavour of plane euclidian geometry. For clarity reasons we shall
describe $B_{0}^{n}(k)$.

Let $\pi_1 :k^n\to k^n$ be the projection of $k^n$ on the $Ox_1$ axis, i.e.
$\pi_1((x_1, x_2, \cdots ,x_n))=(x_1, 0,\cdots,0)$ for all
$(x_1, x_2, \cdots ,x_n)\in k^n$ and $\pi^1:=Id_{k^{n}}-\pi_1$, the
projection of $k^n$ on the hiperplane $x_1=0$, that is
$\pi^{1}((x_1, x_2, \cdots ,x_n))=(0,x_2,\cdots,x_n)$ for all
$(x_1, x_2, \cdots ,x_n)\in k^n$. Then $\pi_1\ot \pi^1$ is a solution
of the Hopf equation and the bialgebra $B_{0}^{n}(k):=B(\pi_1\ot \pi^1)$
can be described as follows:

$\bullet$ $B_{0}^{n}(k)$ is the free algebra generated by
$(c_{ij})_{i,j=1,\cdots, n}$ with the relation
$$
c_{i1}=0, \quad \quad c_{jk}c_{1l}=\delta_{kj}\delta_{l1}c_{11}
$$
for all  $i$, $j\geq 2$ and $k$, $l\geq 1$, where $\delta_{uv}$ is the
Kronecker simbol.

$\bullet$ The comultiplcation $\Delta$ and the counity $\varepsilon$
are given in such a way that the matrix $(c_{ij})_{i,j}$
is comultiplicative.

The proof is similar to the one of proposition \ref{ox}. Among the
elements $(x_{uv}^{ji})$, which define $\pi_1\ot \pi^1$, the only
nonzero elements are
$$
x_{t1}^{t1}=1, \quad \forall t\geq 2.
$$
If $i\neq 1$, all the relations $\chi(i,j,k,l)=0$ are $0=0$,
with the exception of the relations $\chi(i,j,j,1)=0$ for all
$j\geq 2$, which give us $0=c_{i1}$ for all $i\geq 2$.
If $i=1$ the relations
$\chi(1,j,k,l)=0$ give us $c_{jk}c_{1l}=\delta_{kj}\delta_{l1}c_{11}$
for all  $j\geq 2$ and $k$, $l\geq 1$.

New types of bialgebras can be constructed starting from projections of
$k^n$ on different intersections of hyperplanes.
\end{remarks}

\subsection{A five dimensional noncommutative noncocommutative bialgebra}

Let $p$ be a prime number. The classification of $p$ dimensional Hopf
algebras over a field of positive characteristic is still an open problem
(we remind that in \cite{Z} Zhu proved that, over an
algebraically closed field of characteristic zero, any $p$ dimensional
Hopf algebra is isomorphic to the groupal algebra $k[\mathbf{Z}_p]$).
Classifying the $p$ dimensional bialgebras seems to be a much more
complicated problem.

\begin{remark}
We notice that $y$ does not appear in the relations of $E_{0}^{2}(k)$.
As $\Delta(y-1)=(y-1)\ot y +1\ot (y-1)$ and $\varepsilon (y-1)=1$,
we get that the two-sided ideal generated by $y-1$ is also a coideal.
We can add the new relation $y=1$ in the definition of $E_{0}^{2}(k)$
and we obtain a three dimensional noncocommutative
bialgebra. We denote this bialgebra with ${\cal T}(k)$. Then:

$\bullet$ As a vector space, ${\cal T}(k)$ is three dimensional
with $\{1,x,z \}$ a $k$-basis.

$\bullet$ The multiplication rule is given by:
$$
x^2=x, \quad xz=zx=z^2=0.
$$
$\bullet$ The comultiplication $\Delta$ and the counity $\varepsilon$
are given by
$$
\Delta(x)=x\ot x, \quad \Delta(z)=x\ot z+z\ot 1,\quad
\varepsilon(x)=1, \quad \varepsilon(z)=0.
$$
In \cite{Ka}, over a field of characteristic two, two examples of three
dimensional bialgebras are given. Both of them are commutative and
cocommutative. Our ${\cal T}(k)$ differs from one of them only by the
relation $\Delta(z)=x\ot z + z\ot 1$
(respectively $\Delta(z)=1\ot z +z\ot 1$ in \cite{Ka}).
This minor change of $\Delta$ (in our case $k$ being a field of
arbitrary characteristic) makes ${\cal T}(k)$ noncocommutative.

${\cal T}(k)^{*}$ is a three dimensional noncommutative bialgebra. It
follows that ${\cal T}(k)\ot {\cal T}(k)^{*}$ is a nine dimensional
noncommutative and noncocommutative bialgebra.
\end{remark}

If in the preceding remark we have constructed a noncocommutative but
commutative three dimensional bialgebra, now we shall construct a special
$R$ such that $B(R)$ is a five dimensional noncommutative
noncocommutative bialgebra.

\begin{proposition}
Let $k$ be a field and $R$ be the matrix of ${\cal M}_{4}(k)$ given by
$$
R=
\left(
\begin{array}{cccc}
1&0&0&0\\
0&1&1&0\\
0&0&1&0\\
0&0&0&1
\end{array}
\right)
$$
Then:
\begin{enumerate}
\item $R$ is a solution of the Hopf equation if and only if $k$ has the
characteristic two. In this case $R$ is commutative.
\item If $char(k)=2$, then the bialgebra $B(R)$ is the free algebra
generated by $x$, $y$, $z$ with the relations:
$$
x^2=x, \quad y^2=z^2=yx=yz=0, \quad xy=y, \quad xz=zx=z.
$$
The comultiplication $\Delta$ and the counity $\varepsilon$ are given by:
$$
\Delta(x)=x\ot x+y\ot z, \quad
\Delta(y)=x\ot y+y\ot x+ y\ot zy ,\quad
\Delta(z)=z\ot x+x\ot z+zy\ot z,
$$
$$
\varepsilon(x)=1, \quad \varepsilon(y)=\varepsilon(z)=0.
$$
Furthermore, $dim_k(B(R))=5$.
\end{enumerate}
\end{proposition}

\begin{proof}
1. By a direct computation we obtain:
$$
R^{12}=
\left(
\begin{array}{cccccccc}
1&0&0&0&0&0&0&0\\
0&1&0&0&0&0&0&0\\
0&0&1&0&1&0&0&0\\
0&0&0&1&0&1&0&0\\
0&0&0&0&1&0&0&0\\
0&0&0&0&0&1&0&0\\
0&0&0&0&0&0&1&0\\
0&0&0&0&0&0&0&1
\end{array}
\right)
\quad
R^{23}=
\left(
\begin{array}{cccccccc}
1&0&0&0&0&0&0&0\\
0&1&1&0&0&0&0&0\\
0&0&1&0&0&0&0&0\\
0&0&0&1&0&0&0&0\\
0&0&0&0&1&0&0&0\\
0&0&0&0&0&1&1&0\\
0&0&0&0&0&0&1&0\\
0&0&0&0&0&0&0&1
\end{array}
\right)
$$
$$
R^{13}=
\left(
\begin{array}{cccccccc}
1&0&0&0&0&0&0&0\\
0&1&0&0&1&0&0&0\\
0&0&1&0&0&0&0&0\\
0&0&0&1&0&0&1&0\\
0&0&0&0&1&0&0&0\\
0&0&0&0&0&1&0&0\\
0&0&0&0&0&0&1&0\\
0&0&0&0&0&0&0&1
\end{array}
\right)
$$
It follows that:
$$
R^{12}R^{23}=
\left(
\begin{array}{cccccccc}
1&0&0&0&0&0&0&0\\
0&1&1&0&0&0&0&0\\
0&0&1&0&1&0&0&0\\
0&0&0&1&0&1&1&0\\
0&0&0&0&1&0&0&0\\
0&0&0&0&0&1&1&0\\
0&0&0&0&0&0&1&0\\
0&0&0&0&0&0&0&1
\end{array}
\right)
\quad \mbox{and} \quad
R^{23}R^{13}R^{12}=
\left(
\begin{array}{cccccccc}
1&0&0&0&0&0&0&0\\
0&1&1&0&\alpha&0&0&0\\
0&0&1&0&1&0&0&0\\
0&0&0&1&0&1&1&0\\
0&0&0&0&1&0&0&0\\
0&0&0&0&0&1&1&0\\
0&0&0&0&0&0&1&0\\
0&0&0&0&0&0&0&1
\end{array}
\right)
$$
where $\alpha =1+1$. Hence, $R$ is a solution of the Hopf equation if and
only if $\mbox{char}(k)=2$.
In this case we also have that
$$
R^{12}R^{13}=R^{13}R^{12}=
\left(
\begin{array}{cccccccc}
1&0&0&0&0&0&0&0\\
0&1&0&0&1&0&0&0\\
0&0&1&0&1&0&0&0\\
0&0&0&1&0&1&1&0\\
0&0&0&0&1&0&0&0\\
0&0&0&0&0&1&0&0\\
0&0&0&0&0&0&1&0\\
0&0&0&0&0&0&0&1
\end{array}
\right)
$$
i.e. $R$ is commutative.

2. Suppose that $\mbox{char}(k)=2$. Starting as in the proof of the
above propositions, with a basis in a two dimensional vector space, we
obtain that among $(x_{ij}^{kl})$ the only nonzero elements
are
$$
x_{11}^{11}=x_{22}^{22}=x_{21}^{21}=x_{12}^{12}=x_{12}^{21}=1.
$$
Now, the relations $\chi(i,j,k,l)=0$, written in the lexicografical
order are:
$$
c_{11}c_{11}=c_{11}, \quad c_{11}c_{12}=c_{12}, \quad
c_{12}c_{11}=0, \quad  c_{12}c_{12}=0,
$$
$$
c_{21}c_{11}+c_{11}c_{21}=0, \quad c_{21}c_{12}+c_{11}c_{22}=c_{11},
$$
$$
c_{22}c_{11}+c_{12}c_{21}=c_{11}, \quad c_{22}c_{12}+c_{12}c_{22}=c_{12},
$$
$$
c_{11}c_{21}=c_{21}, \quad c_{11}c_{22}=c_{22}, \quad
c_{12}c_{21}=0, \quad c_{12}c_{22}=0,
$$
$$
c_{21}c_{21}=0, \quad c_{21}c_{22}=c_{21}, \quad
c_{22}c_{21}=c_{21}, \quad c_{22}c_{22}=c_{22}.
$$
Now, if we denote $c_{11}=x$, $c_{12}=y$, $c_{21}=z$, $c_{22}=t$ and
using that $\mbox{char}(k)=2$ we get the following relations:
$$
x^2=x, \quad y^2=z^2=yx=yz=0, \quad xy=y, \quad xz=zx=z,
$$
$$
zy=x+t, \quad t^2=t, \quad xt=t, \quad tx=x,
$$
$$
yt=0, \quad ty=y, \quad zt=tz=z.
$$
So, $t$ is in the free algebra generated by $x$, $y$, $z$ and
$$
t=zy-x.
$$
If we substitute $t$ in all the relations in which $t$ is involved,
then these become identities. The relations given in the
statement of the proposition remain. The formula for $\Delta$ follows, as
the matrix
$$
\left(
\begin{array}{cc}
x&y\\
z&t
\end{array}
\right)
$$
was comultiplicative.

We shall prove now that $\mbox{dim}_k(B(R))=5$, more
exactly we will show in an elementary way (without the diamond lemma) that
$\{1, x, y, z, zy \}$ is a $k$-basis for $B(R)$. From the relations which
define $B(R)$ we obtain:
$$
x(zy)=zy, \quad
(zy)^2=(zy)x=y(zy)=(zy)y=z(zy)=(zy)z=0.
$$
All these relations give us that $\{1, x, y, z, zy \}$ is a sistem
of generators of $B(R)$ as a vector space over $k$. It remain to check that
they are liniar independent over $k$. Let $a$, $b$, $c$, $d$, $e\in k$
such that
$$
a+bx+cy+dz+e(zy)=0.
$$
First, we multiply to the left with $y$ and we get $a=0$. Then we multiply
to the right with $z$ and we obtain that $b=0$. Now we multiply
to the right with $x$ and we get $d=0$. If we multiply now to the left
with $z$ we get $c=0$, and $e=0$ follows.
Hence $\{1, x, y, z, zy \}$ is a $k$-basis for $B(R)$.
\end{proof}

\begin{remarks}
1. The bialgebra $B(R)$ constructed in the above proposition is not a Hopf
algebra. We can localize $B(R)$ and we get a Hopf algebra which is
isomorphic to the grupal Hopf algebra $k[G]$, where $G$ is a group with
two elements. Indeed, let $S$ be a potential antipode. Then:
$$
S(x)x+S(y)z=1, \quad S(x)y+S(y)t=0.
$$
If we multiply the second equation to the right with $z$ we get
$S(y)z=0$, so $S(x)x=1$. But $x^2=x$, so $x=1$ and then $y=0$, $t=1$.
We obtain the Hopf algebra $k<z\mid z^2=0>$,
$\Delta(z)=z\ot 1+1\ot z$, $\varepsilon(z)=0$. If we denote
$g=z+1$ then $g^2=1$, $\Delta(g)=g\ot g$, $\varepsilon(g)=1$, hence
$B(R)$ is isomorphic to the Hopf algebra $k[G]$, where $G=\{1, g\}$
is a group with two elements.

2. Directly from the proof, we obtain an elementary definition for
the bialgebra $B(R)$:

$\bullet$ As a vector space $B(R)$ is five dimensional with
$\{1,x,y,z,t \}$ a $k$-basis.

$\bullet$ The multiplication rule is given by:
$$
x^2=x, \quad y^2=z^2=0, \quad t^2=t,
$$
$$
xy=y,\quad yx=0,\quad xz=zx=z,\quad xt=t,\quad tx=x,
$$
$$
yz=0,\quad zy=x+t,\quad yt=0,\quad ty=y, \quad zt=tz=z.
$$
$\bullet$ The comultiplcation $\Delta$ and the counity $\varepsilon$
are given in such way that the matrix
$$
\left(
\begin{array}{cc}
x&y\\
z&t
\end{array}
\right)
$$
is comultiplicative.

3. As $R$ is commutative we can construct the bialgebra $\overline{B}(R)$:
it is the quotient  $k[X, Z]/(X^2-X, Z^2, XZ-Z)$ of the polinomial bialgebra
$k[X, Z]$ with the coalgebra structure given by
$$
\Delta(X)=X\ot X,\quad \Delta(Z)=X\ot Z + Z\ot X
$$
$$
\varepsilon(X)=1,\quad \varepsilon(Z)=0.
$$

4. Recently, in \cite{IS}, considering quotients of the bialgebras $B(R)$
for various solutions $R$ of the Hopf-equation, numerous examples of finite
dimensional noncommutative and noncocommutative bialgebras are constructed.
We shall present one of them, which is a quotient of $E_{0}^2(k)$.

Let $n\geq 2$ be a natural number. Then, the two-sided ideal $I$ of
$E_{0}^{2}(k)$ generated by $y^n-y$, $zy$, $xy-x$ and $yx-x$ is a biideal
and
$B_{2n+1}(k):=E_{0}^{2}(k)/I$ is a $2n+1$-dimensional noncommutative
noncocommutative bialgebra (see \cite{IS} for the proof). The bialgebra
$B_{2n+1}(k)$ can be described as follows:

$\bullet$ $B_{2n+1}(k)$ is the free algebra generated by $x$, $y$, $z$ with
the relations:
$$
x^2=x,\quad xz=zx=z^2=zy=0,\quad y^n=y,\quad xy=yx=x.
$$
$\bullet$ The comultiplication $\Delta$ and the counity $\varepsilon$
are given by:
$$
\Delta(x)=x\ot x,\quad \Delta(y)=y\ot y,\quad \Delta(z)=x\ot z+z\ot y
$$
$$
\varepsilon(x)=\varepsilon(y)=1,\quad \varepsilon(z)=0.
$$
\end{remarks}

{\sl Acknowledgement.} The author thanks Bogdan Ion and Mona Stanciulescu for
their comments on the preliminary version of this paper.
We also want to thank the referee for his/her valuable suggestions.

\end{document}